\documentclass[12pt]{amsart}
\usepackage{amsmath,amssymb,amsfonts,amsthm,amsopn}
\usepackage{graphics}

 \def\Spnr{Sp(d,\R)}
 \def\Gltwonr{GL(2d,\R)}

\newcommand{\tfa}{time-frequency analysis}

\newcommand{\ft}{Fourier transform}
\newcommand{\stft}{short-time Fourier transform}

\newcommand{\tf}{time-frequency}

\newcommand{\fif}{if and only if}
\newcommand{\tfs}{time-frequency shift}

\newcommand{\modsp}{modulation space}
\newcommand{\psdo}{pseudodifferential operator}

\newtheorem{theorem}{Theorem}[section]
\newtheorem{lemma}[theorem]{Lemma}
\newtheorem{corollary}[theorem]{Corollary}
\newtheorem{proposition}[theorem]{Proposition}
\newtheorem{definition}[theorem]{Definition}
\newtheorem{cor}[theorem]{Corollary}
\newtheorem{remark}[theorem]{Remark}

\newcommand{\beqa}{\begin{eqnarray*}}
\newcommand{\eeqa}{\end{eqnarray*}}

\DeclareMathOperator*{\supp}{supp}

\newcommand{\field}[1]{\mathbb{#1}}
\newcommand{\bR}{\field{R}}        
\newcommand{\bN}{\field{N}}        
\newcommand{\bZ}{\field{Z}}        
\def\G{\mathcal{G}}

\def\la{\lambda}

\def\cS{\mathcal{S}}

\def\cG{\mathcal{G}}
\def\cM{\mathcal{M}}

\def\cC{\mathcal{C}}

\def\a{\aleph}

\def\rd{\bR^d}

\def\rdd{{\bR^{2d}}}
\def\zdd{{\bZ^{2d}}}

\def\lrd{L^2(\rd)}

\def\zd{\bZ^d}

\def\intrd{\int_{\rd}}
\def\intrdd{\int_{\rdd}}

\def\R{\right)}
\def\l{\langle}
\def\r{\rangle}
\def\<{\left<}
\def\>{\right>}

\def\inv{^{-1}}

\def\mv1{M_v^1}

\def\phas{(x,\o )}
\def\mn{(m,n)}
\def\mn'{(m',n')}

\def\Spnr{Sp(d,\R)}

\def\o{\eta}
\def\a{\alpha}
\def\b{\beta}

\def\R{\mathbb{R}}
\def\Ren{\mathbb{R}^d}

\def\f{\varphi}

\def\Sn2{S_{2}(L^{2}(\Ren))}
\def\S1{S_{1}(L^{2}(\Ren))}
\def\sig00{\sigma_{0,0}}

\def\la{\langle}
\def\ra{\rangle}

\begin{document}

\begin{abstract}
A general principle says that the matrix of a Fourier integral
operator with respect to wave packets is concentrated near the curve
of propagation. We prove a precise version of this principle for
Fourier integral operators  with a smooth phase and a symbol in the
Sj\"ostrand class and use Gabor frames as wave packets. The almost
diagonalization of such Fourier integral operators suggests a specific
approximation by  (a sum of) elementary operators, namely   modified Gabor multipliers.
We derive error estimates for such approximations. The methods are
taken from \tfa  .
\end{abstract}

\title{Approximation of Fourier Integral
Operators by Gabor multipliers}

\author{Elena Cordero, Karlheinz Gr\"ochenig and Fabio Nicola}
\address{Department of Mathematics,
University of Torino, via Carlo Alberto 10, 10123 Torino, Italy}
\address{Faculty of Mathematics,
University of Vienna, Nordbergstrasse 15, A-1090 Vienna, Austria}
\address{Dipartimento di Matematica,
Politecnico di Torino, corso Duca degli Abruzzi 24, 10129 Torino,
Italy}

\email{elena.cordero@unito.it}
\email{karlheinz.groechenig@univie.ac.at}
\email{fabio.nicola@polito.it}
\thanks{K.\ G.\ was
  supported in part by the  project P22746-N13  of the
Austrian Science Foundation (FWF)}
\maketitle

\newcommand{\fio}{Fourier integral operator}

\section{Introduction}

A fundamental  principle expressed  by Cordoba and Fefferman~\cite{CF78}
says that   Fourier integral operators map wave packets to wave
packets. Furthermore each wave packet is  transported 
according to the canonical flow in phase space that is  associated to the operator.  

Rigorous versions of this principle were proved for various types of
\fio s and wave packets. 
Cordoba and Fefferman used generalized Gaussians as wave packets and
considered \fio s with constant symbol. In the last decade the
transport of wave packets by \fio s was investigated for
  curvelet-like frames by Smith~\cite{smi98}, for curvelets by Candes
  and Demanet~\cite{CD05},   for
shearlets by Guo and Labate~\cite{GL08}, and  more recently for Gabor
frames by Rodino and two us~\cite{fio1}. 

In this paper we investigate \fio s with a ``tame'' phase and
non-smooth symbols with respect to \tfs s and Gabor
frames. Following~\cite{AS78}, 
we will assume that the  phase function $\Phi $ on $\rdd $ satisfies  the conditions 
 (i) $\Phi\in \cC^{\infty}(\rdd)$,  (ii) 
$\partial_z^\a \Phi \in L^\infty (\rdd )$ for $\a \geq 2$, and   (iii)
$\inf _{x,\eta \in \rdd}    |\det\,\partial^2_{x,\eta}
\Phi(x,\o)|\geq \delta >0$. For brevity, we call $\Phi$  a tame phase
function.   Then the \fio\ $T$ with phase function $\Phi $ and the  symbol 
$\sigma $  is formally defined to be 
\begin{equation}
  \label{eq:i1}
   Tf(x)=\intrd e^{2\pi i \Phi\phas} \sigma\phas \hat{f}(\o)d\o.
\end{equation}
\fio s with a tame phase were investigated by many authors, both in
hard analysis~\cite{Bou97} and in \tfa  ~\cite{CT09,fio1,fio2}. They arise  the study of Schr\"odinger
operators, for instance, for the description of the resolvent  of the
Cauchy problem for the Schr\"odinger equation with a quadratic type
Hamiltonian~\cite{AS78, helffer84, helffer-rob1}. In our preceeding work~\cite{fio1,fio2,fio5}  we have adopted
the point of view that  \fio s with a tame phase  are best studied
with those \tf\ methods that correspond to a constant geometry of the
wavepackets. Precisely,  let $z= (x,\eta ) \in \rdd $ be a point in phase space
and 
$$
\pi (z) f(t) = e^{2\pi i \eta \cdot t} f(t-x), \qquad t\in \rd $$
be the corresponding phase space shift (or \tfs ).
Roughly speaking, these wave packets correspond to a uniform partition of the
phase space.

In this paper we continue the investigation of the matrix of a \fio\
with respect to (frames of)   phase space shifts. 
A  first version
of a \tfa\ of \fio s was proved in ~\cite{fio1}.  The formulation is
in the spirit of Candes and Demanet and is intriguing in its
simplicity. 

\begin{theorem} \label{old}
   Assume that $\Phi $ is a tame phase function and $\chi $ is the
corresponding     canonical transformation of $\rdd $. If    $g\in \cS (\rd
 )$ and if $\sigma $ belongs
 to the Sobolev space $ W^\infty
 _{2N}(\rdd)$,  
 then there exists a  constant
$C_N>0$ such that
\begin{equation}\label{5f}
   |\langle T \pi(\mu)g, \pi(\lambda)g\rangle|\leq C_N {\la \chi(\mu)
   -\lambda\ra^{-2N}
   },\qquad \forall \lambda,\mu \    \in\R^{2d} \, .
\end{equation}
\end{theorem}

Our first goal is to sharpen Theorem~\ref{old} in several respects.

(a) We keep the assumptions on the phase function $\Phi $ (tame 
phase), but allow non-smooth symbols taken in a generalized Sj\"ostrand
class. 

(b) We enlarge the class of ``basis'' functions or ``windows'' $g$  and extend
Theorem~\ref{old} to hold for windows in certain \modsp s $M^1_v(\rd
)$. 

 The appearance of \modsp s is not longer a surprise.    It is  well
 established   that \modsp s are the appropriate
function spaces for phase space analysis and \tfa\ with constant
geometry. In particular  the Sj\"ostrand class and
its generalizations  are   a special case of the
 \modsp s and they  arise as suitable  symbol classes for \psdo s. 
The use of the Sj\"ostrand class for \fio s goes
back to Boulkhemair~\cite{Bou97} and features prominently in 
recent work on \psdo s and \fio s (see~\cite{Gro06ped} for a survey). 
Furthermore, the  theory of \tfa\ and  \modsp s is now well developed
and offers  new tools  for the investigation of \fio s, in addition
to the classical methods.

The arguments and  proofs of  the extension of Theorem~\ref{old} are
in the spirit of \cite{fio1,gro06}, but require a substantial amount
of technicalities. A main technical result is a new characterization of
  the Sj\"ostrand class in Proposition~\ref{proposition4}. This
  characterization 
is perfectly adapted  to the  investigation of \fio s and measures  \modsp\ norms with a
variable family of windows rather than a fixed window. 

Our second goal is to derive  a discretized version of
Theorem~\ref{old}. We will show that  a
\fio\ can be approximated  by (sums of) elementary operators of the
form 
\begin{equation}
  \label{eq:i2}
  M_{\bf a} f=M_{\bf
    a}^{\chi',g,\Lambda}f=\sum_{\lambda\in\Lambda} a_{\lambda} \langle
  f,\pi(\lambda) g\rangle\pi(\chi'(\lambda))g \, ,
\end{equation}
where $\Lambda $ is a lattice in $\rdd$, $(a_{\lambda})$ a symbol
sequence on $\Lambda $,  and $\chi ' :\Lambda \to
\Lambda $ is  a discrete
version of the canonical transformation associated to a phase function
$\Phi $.  
Such operators are suggested by a recent approximation theory for
\psdo s through Gabor multipliers~\cite{approx}. The modified Gabor
multipliers in~\eqref{eq:i2} are  adapted precisely to the canonical
transformation.  We will approximate a \fio\ with tame phase by
modified Gabor multipliers and  prove an 
approximation theorem for \fio s (Theorem~\ref{teorema2.1}). 

The approximation theory also yields an alternative proof for the boundedness of
\fio s on \modsp s. The boundedness and Schatten class properties of \fio s were studied
 under various assumptions on the phase and the symbol. We refer to
 the book~\cite{stein93} and the
 articles~\cite{Bou97,CF78,CT09,fio5,RS06} for a sample of  contributions.

 Gabor multipliers have many applications in signal
processing and acoustics~\cite{daub88,torr10} and are especially
useful for the numerical realization of \psdo s. For the study of \fio s
we introduce  a
new type of Gabor multipliers. We hope that the modified Gabor
multipliers will also prove useful for the numerical realization and
approximation of \fio s.

Our paper is organized as follows: In Section~2 we collect some
concepts and the necessary facts  from \tfa . Section~3 is devoted to the
almost diagonalization of \fio s with a tame phase. We will prove a
substantial generalization  of Theorem~\ref{old}. In Section~4 we define
and study the modified Gabor multipliers mentioned in~\eqref{eq:i2} and
prove their boundedness on \modsp s. In Section~5 we study the
representation and approximation of \fio s by sums of modified Gabor
multipliers and derive a quantitative approximation theorem. Further
remarks are contained in Section~6.

\section{Notation and preliminary results}
We recall some notation and tools from time-frequency analysis. For an
exposition and details we refer to the book~\cite{book}. \par

(i) First we define the translation and modulation operators
$$ T_xf(t)=f(t-x)\quad{\rm and}\quad M_{\o}f(t)= e^{2\pi i \o  t}f(t),$$
for $\lambda=(x,\o)\in\rdd$, $f\in\cS'(\rd)$. Their composition is the
\tfs\
$$\pi(\lambda)= M_{\o}T_x.
$$
(ii) Now  consider a distribution $f\in\cS '(\rd)$ and a Schwartz function $g\in\cS(\rd)\setminus\{0\}$ (which will be called {\it window}).
The short time Fourier transform of $f$ with respect to $g$ is given by
\[
V_g f(x,\o)=\langle f,M_{\o}T_x g\rangle=\langle f,\pi(\lambda)g
\rangle= \int _{\rd} f(t)\overline{g(t-x)}e^{-2\pi it\cdot\o}\,dt,
\quad {\rm for}\ x,\eta \in\R^d \, .
\]
The  last integral makes sense for dual pairs of function spaces,
e.g., for  $f,g\in
L^2(\rd)$.  The bracket $\langle \cdot , \cdot \rangle$ makes sense for
dual pairs of function or distribution spaces, in particular for $f\in
\cS ' (\rd )$ and $g\in \cS (\rd )$. 
\par

(iii)  In the sequel, we set $v\phas=v_s\phas=\la
\phas\ra^s=(1+|x|^2+|\o|^2)^{s/2}$ and we denote by
$\mathcal{M}_v(\rdd)$ the space of $v$-moderate weights on $\rdd$;
these  are measurable functions $m>0$ satisfying $m(z+\zeta)\leq C
v(z)m(\zeta)$ for every $z,\zeta\in\rd$.\par

(iv) Then, for $1\leq p \leq \infty$ and $m\in
\mathcal{M}_v(\rdd)$ we denote by $M ^{p}_m(\R^d)$ the space of
distributions $f\in\cS'(\rd)$ such that
\[
\|f\|_{M ^{p}_m(\R^d)}:=\|V_g f\|_{L^p_m(\R^{2d})}<\infty.
\]
This definition is meaningful and does not depend on the choice of
the window $g\in \cS (\rd ), g \neq 0$. If fact, different $g\in
M^1_v(\rd )$ yield equivalent norms on $M^p_m (\rd )$ whenever
$m\in \cM _v(\rdd)$.

(v) Modulation spaces possess a  discrete description as well.  Consider a lattice of the form
$\Lambda=A\zdd$  with $A\in GL(2d,\R)$. The collection of
time-frequency shifts $\G(g,\Lambda)=\{\pi(\lambda)g:\
\lambda\in\Lambda\}$ for a  non-zero $g\in L^2(\rd)$ is called a
Gabor system. The set $\G(g,\Lambda)$ is a Gabor frame, if there exist 
constants $A,B>0$ such that
\begin{equation}\label{gaborframe}
A\|f\|_2^2\leq\sum_{\lambda\in\Lambda}|\langle f,\pi(\lambda)g\rangle|^2\leq B\|f\|^2_2\qquad \forall f\in L^2(\rd).
\end{equation}
Gabor frames give the following characterization of the Schwartz space
$\cS(\rd)$  and of the modulation spaces $M^{p}_m(\rd)$.  Fix $g\in \cS (\rd ), g \neq 0$, then
\begin{align}
 \label{osservazione0}
f\in\cS(\rd) & \Leftrightarrow \sup_{\lambda\in\Lambda}(1+|\lambda|)^N
\langle f,\pi(\lambda)g\rangle<\infty \quad\forall N\in\bN,\\
\label{osservazione0bis}
f\in M^{p}_m(\rd) &\Leftrightarrow \Big(\sum_{\lambda\in\Lambda}|\langle f,\pi(\lambda)g\rangle|^p m(\lambda)^p\Big)^{1/p}<\infty,
\end{align}
and the latter expression gives an equivalent norm for
$M^p_m(\rd)$. 

\par Finally we
recall that, if $A=B=1$ in \eqref{gaborframe}, then
$\G(g,\Lambda)$ is called a Parseval frame, and \eqref{gaborframe}
implies the expansion
\begin{equation}\label{parsevalframe}
f=\sum_{\lambda\in\Lambda}\langle f,\pi(\lambda)g\rangle\pi(\lambda)g\qquad\forall f\in L^2(\rd)
\end{equation}
with unconditional convergence in $L^2(\rd)$.  If $f\in M^p_m(\rd
)$ and $g\in M^1_v(\rd )$ ($m\in\mathcal{M}_v(\rdd)$), then this
expansion converges unconditionally  in $M^p_m(\rd )$ for $1\leq
p<\infty$.

(vi) Amalgam spaces:  Let $m\in \cM _v(\rdd)$. The amalgam space
$W(L^p)(\rd ) $ consists of all essentially bounded functions such
that the norm
$$
\|f\|_{W(L^p_m)} = \sum _{k\in \zd } \sup _{x\in [0,1]^d} |f(x+k)|
m(k) < \infty
$$
is finite.  If $f\in W(L^p_m)$ is also continuous and $\Lambda $ is a
lattice in $\rd$, then the sequence
$(f(\lambda ))_{\lambda \in \Lambda }$ is in the sequence space $\ell
^p_m(\Lambda )$.
We will use several times the usual convolution relations for amalgam
spaces
\begin{equation}
 \label{eq:may2}
L^p_m \ast W(L^1_v) \subseteq W(L^p_m)   \, ,
\end{equation}
which hold for $1\leq p\leq \infty$ and all  $m\in \cM _v(\rdd)$.

\section{Almost diagonalization of FIOs}\label{ad}
For a given function $f$ on $\rd$,  the
Fourier integral operator (FIO in short) with symbol $\sigma$ and phase
$\Phi$ can be formally defined by
\begin{equation}\label{fio}
   Tf(x)=\intrd e^{2\pi i \Phi\phas} \sigma\phas \hat{f}(\o)d\o.
\end{equation}
To avoid technicalities we take $f\in\mathcal{S}(\R^d)$ or, more
generally, $f\in M^1(\rd )$. If $\sigma\in L^\infty(\rdd )$ and if
the phase $\Phi$ is real, the integral converges absolutely and
defines a function in $L^\infty(\rd )$.\par We will  consider  a
class of Fourier integral operators which arises in the study of
Schr\"odinger type equations. This class of operators was already
studied in the spirit of time-frequency analysis in \cite{fio1}.
\begin{definition}
A phase
function $\Phi(x,\eta)$ is called \emph{tame}, if $\Phi $  satisfies  the following
properties:\par\medskip
\noindent
{\bf (i)} $\Phi\in \cC^{\infty}(\rdd)$;\\
{\bf (ii)} for $z=\phas$,
\begin{equation}\label{phasedecay}
|\partial_z^\a \Phi(z)|\leq
C_\a,\quad |\a|\geq
2;\end{equation}
{\bf  (iii)}
there exists $\delta>0$ such
that
\begin{equation}\label{detcond}
   |\det\,\partial^2_{x,\eta} \Phi(x,\o)|\geq \delta.
\end{equation}
\end{definition}
\par
If we set
\begin{equation}\label{cantra} \left\{
                \begin{array}{l}
                y=\nabla_{\eta}\Phi(x,\eta)
                \\
               \xi=\nabla_{x}\Phi(x,\eta), \rule{0mm}{0.55cm}
                \end{array}
                \right.
\end{equation}
we can solve with respect to $(x,\xi)$ by the global inverse
function theorem (see e.g. \cite{krantz}), obtaining a mapping $\chi$
defined by $(x,\xi)=\chi(y,\o)$. As observed in \cite{fio1}, $\chi $
is a smooth bi-Lipschitz canonical
transformation. This means that\par\medskip\noindent
-- $\chi$ is a  smooth diffeomorphism on $\rdd$;\\
-- both $\chi$ and $\chi^{-1}$ are Lipschitz continuous;\\
-- $\chi$ preserves the symplectic form. That is, if $(x,\xi)=\chi(y,\eta)$,

\begin{equation}\label{simp}
\sum_{j=1}^d dx_j\wedge d\xi_j= \sum_{j=1}^d dy_j\wedge d\eta_j.
\end{equation}
Observe that condition \eqref{simp} is equivalent to saying that the differential $d\chi(y,\eta)$, at every point
 $(y,\eta)\in\R^{2d}$, as a linear map $\R^{2d}\to\R^{2d}$, is represented by a symplectic matrix, i.e. belonging to the group $$\Spnr:=\left\{A\in\Gltwonr:\;^t\!AJA=J\right\},$$  where
$$
J=\bmatrix 0&I_d\\-I_d&0\endbmatrix.
$$
In particular, when $\chi$ is linear, this means that $\chi$ itself is represented by a symplectic matrix. \par


We now address  to the
problem of the almost diagonalization of Fourier integral
operators with a tame phase function. It turns out that, when the symbol
is regular enough, the matrix representation of such an operator
is almost diagonal with respect to a Gabor system.\par


 Theorem~\ref{old}  in  the introduction  was the first  precise
 result about the almost diagonalization of \fio s in \tfa :
\emph{ Assume that $\Phi $ is a tame phase function and   $g\in \cS (\rd
 )$.  If $\sigma $ belongs
 to the Sobolev space $ W^\infty
 _{2N}(\rdd)$ (consisting of all distributions with essentially bounded
 derivatives up to order $2N$), then 
there exists a  constant
$C_N>0$ such that}
\begin{equation}\label{5fa}
   |\la T \pi(\mu)g, \pi(\lambda)g\ra|\leq C_N {\la \chi(\mu)
   -\lambda\ra^{-2N}
   },\qquad \forall \lambda,\mu \    \in\R^{2d},
\end{equation}
\emph{where $\chi$ is the canonical
transformation generated by $\Phi$.}


This statement  proved in~\cite[Theorem 1]{fio1}. For the extension of
Theorem~\ref{old} we next introduce the appropriate symbol class.   
If $m\in\mathcal{M}_{v_{s}}(\R^{2d})$, with $v_{s}=v_{s}(x,\o)$,
$(x,\o)\in\rdd$, $s\in\R$, we consider the symbol class 
$M^{\infty,1}_{1\otimes m}(\R^{2d})$ of $\sigma\in\cS'(\R^{2d})$
satisfying
\[
\|\sigma\|_{M^{\infty,1}_{1\otimes m}(\R^{2d})}=\|V_{\Psi}
\sigma\|_{L^{\infty,1}_{1\otimes m}(\R^{4d})}=\int_{\rdd}
\sup_{z\in\rdd}|V_{\Psi}(z,\zeta)|m(\zeta)\,d\zeta<\infty.
\]
In \tfa\ this symbol class  is just  a special case of a \modsp\ ~\cite{F1}, in theory
of \psdo s $M^{\infty,1}_{1\otimes m}$ is often referred to as a
generalized Sj\"ostrand class after ~\cite{Sjo94}. 

For FIOs with symbols in $M^{\infty,1}_{1\otimes m}(\R^{2d})$
there is an  almost diagonalization result  similar  to
Theorem \ref{old}.\par
The window $g$ will be chosen sufficiently regular, as follows. Given $s\geq0$, let
\begin{equation}\label{g1}
N\in\mathbb{N},\quad N>\frac{s}{2}+d,\quad{\rm and}\quad
  g\in M^1_{v_{{4N}}\otimes v_{{4N}}}(\rd),
  \end{equation}
here $v_{4N}$ is  a weight function defined on $\rd$.

\begin{theorem}\label{6bis}
Let $m\in\mathcal{M}_{v_{s}}(\rdd)$, $s\geq0$, $N$, $g$ satisfying \eqref{g1} and
such that $\G(g,\Lambda)$ is a Parseval frame for $L^2(\rd)$. If
the phase $\Phi$ satisfies  conditions {\bf (i)}, {\bf (ii)} and $\sigma\in
M^{\infty,1}_{1\otimes m}(\R^{2d})$
then there exists $H\in W(L^{1}_m)(\R^{2d})$ such that
\begin{equation}\label{boas}
|\langle T \pi(x,\o)g,\pi(x',\o')g\rangle|\leq
H(\o'-\nabla_x\Phi(x',\o),x-\nabla_\eta\Phi(x',\o)),
\end{equation}
for every $ \phas,(x',\o')\in\rdd$, with $\|H\|_{
W(L^{1}_m)}\lesssim \|\sigma\|_{M^{\infty,1}_{1\otimes m}}$
\end{theorem}

Observe that Theorem \ref{6bis} says we can control the
coefficient matrix \linebreak $|\langle T\pi(x,\o)g,\pi(x',\o')g \rangle|$ by
a function $H$ of the difference
$(\o'-\nabla_x\Phi(x',\o),x-\nabla_\eta\Phi(x',\o))$.  It will be
useful in the sequel to have a control of the matrix which depends
on the difference $(\o'-\chi_2\phas,x'-\chi_1\phas)$. This is
achieved  by adding the assumption ${\bf (iii)}$ of ~\cite{AS78} on the phase and
by choosing  another, even larger  modulation space as the  symbol class.

\begin{theorem}\label{6biss}
Let $s\geq0$, $N$, $g$ satisfying \eqref{g1} and
assume  that $\G(g,\Lambda)$ is a Parseval frame for $L^2(\rd)$. If
the phase $\Phi$ is tame
and $\sigma\in M^{\infty}_{1\otimes v_s}(\R^{2d})$,
then there exists $C>0$ such that
\begin{equation}\label{bo}
|\langle T \pi(\mu )g,\pi(\lambda )g\rangle|\leq C \frac{
\|\sigma\|_{M^{\infty}_{1\otimes v_s}}} {\langle \chi (\mu) - \lambda
\rangle ^s},\quad\forall \lambda , \mu \in\rdd.
\end{equation}
\end{theorem}

Theorem~\ref{6biss} improves Theorem~\ref{old} in several
respects: If $s=2N$, then $W^\infty _{2N}(\rdd ) \subseteq
M^\infty _{1\otimes
  v_{2N}}(\rdd )$ \cite[Prop. 6.7]{F1} and $\cS (\rd ) \subseteq  M^1_{v_{{4N}}\otimes
  v_{{4N}}}(\rd)$. Thus we obtain  the same quality of off-diagonal
decay for  significantly larger  classes of symbols and windows.


The remainder of this section is  devoted to set up the tools for the
proof of these theorems.

 For $z,w\in\rdd$, let  $\Phi _{2,z}$  the remainder  in the second
 order Taylor expansion  of the phase
 $\Phi $, i.e.,
\[
\Phi_{2,z}(w)=2\sum_{|\alpha|=2}\int_0^1
(1-t)\partial^\alpha\Phi(z+tw)dt \frac{w^\alpha}{\alpha!} \, .
\]
For  a given window $g$, we set
\begin{equation}\label{psi}
\Psi_z(w)=e^{2\pi i \Phi_{2,z}(w)}
\overline{g}\otimes\widehat{g}(w).
\end{equation}

The main technical work is to show that the set of windows $\Psi _z$
possesses a joint time-frequency envelope. This property will allow us
to replace the $z$-dependent family of windows $\Psi _z$ by a single
window in many estimates.


Before proving the existence of a \tf\  envelope in Lemma
\ref{wienerprop} below, 
we first look at the phase factor
$e^{2\pi i \Psi_{2,z}}$ occurring  in
\eqref{psi}.
\begin{lemma}\label{window} For every $s\in\R$, $N\in\mathbb{N}$,
  $N>\frac{s}{2}+d$,  and $\Psi\in\cS(\rdd)$, we have
\begin{equation}\label{l2}
\sup_{z\in\rdd}|V_\Psi e^{2\pi i \Phi_{2,z}}|\in
L^{\infty,1}_{v_{-4N}\otimes v_{s}}(\R^{4d}),
\end{equation}
with  $v_{-4N}$ and $v_{s}$ being weight functions on
$\rdd$.
\end{lemma}
\begin{proof}
We proceed as in~\cite{fio1}.
Let $\Psi\in\cS(\rdd)$, then
$$|V_\Psi e^{2\pi i \Phi_{2,z}}(u,w)|=\left| \int_\rdd e^{2\pi i
\Phi_{2,z}(\zeta)} T_u\bar{\Psi} (\zeta) e^{-2\pi i \zeta \cdot w}
d\zeta\right|.$$
 Using the identity
$$(1-\Delta_\zeta)^N e^{-2\pi i \zeta w}=\la
   2\pi w\ra^{2N} e^{-2\pi i \zeta w},
$$
we integrate by parts and obtain
$$
 |V_\Psi e^{2\pi i \Phi_{2,z}}(u,w)|\\
   =\frac1{\la
   2\pi w\ra^{2N}}\left|\intrdd (1-\Delta_\zeta)^N\big(e^{2\pi i
\Phi_{2,z}(\zeta)} T_u\bar{\Psi} (\zeta)\big) e^{-2\pi i \zeta \cdot  w}
d\zeta\right|.
$$
By means of Leibniz's formula the factor $(1-\Delta_\zeta)^N
(e^{2\pi i \Phi_{2,z}(\zeta)} T_u\bar{\Psi} (\zeta))$ can be
expressed further as $$e^{2\pi i\Phi_{2,z}(\zeta)}\sum_{|\a|+|\b|\leq 2N}
    p_\alpha(\partial  \Phi_{2,z}(\zeta))
   (T_u\partial_\zeta^\b\bar{\Psi})(\zeta),$$
   where $p_\alpha(\partial \Phi_{2,z}(\zeta))$ is a   polynomial    of derivatives of
   $\Phi_{2,z}$ of degree at most $|\a|$.\par
   As a consequence of~\eqref{phasedecay} we have
$|p_\alpha(\partial \Phi_{2,z}(\zeta))|\leq C_{\a} \langle
\zeta\rangle^{2|\a|}$  for every $z\in\rdd$ with a constant
independent of $z$.  Moreover, the assumption $\Psi \in \cS(\rdd)$
yields that  $$\sup_{|\b|\leq
2N}|T_u\partial_\zeta^\b\bar{\Psi}|\leq C_{N,s} \la
\zeta-u\ra^{-\ell}$$
for every $\ell\geq 0$. Consequently,
\begin{align*}|V_\Psi e^{2\pi i
\Phi_{2,z}}(u,w)|&\lesssim \frac1{\la
   2\pi w\ra^{2N}} \intrdd \sum_{|\a|+|\b|\leq 2N}
     \la \zeta\ra^{2|\a|} \la \zeta
    -u\ra^{-\ell}d\zeta\\
    &\lesssim \frac1{\la
   2\pi w\ra^{2N}}\intrdd \la \zeta\ra^{4N} \la \zeta
    -u\ra^{-\ell}d\zeta\\
    &\lesssim \frac1{\la
   2\pi w\ra^{2N}}\la u\ra^{4N},
\end{align*}
whenever   $\ell> 4N + 2d$. Since $-2N+s<-2d$ by assumption, we obtain
$$\intrdd\sup_{u\in\rdd}\sup_{z\in\rdd}|V_\Psi e^{2\pi i \Phi_{2,z}}|(u,w)\la
u\ra^{-4N}\la w\ra^{s} \,dw<\infty,$$ whence \eqref{l2} follows.
\end{proof}
\begin{remark}\label{osserv}  Notice that
 the two weights $v_{{-4N}}$ and $v_{s}$ compensate each other in \eqref{l2}.
 It is easy to see that some form of compensation is necessary.  For example,
if  $\Phi(\zeta)=|\zeta|^{2}/2$, then
$\Phi_{2,z}(\zeta)=|\zeta|^{2}/2$ independent of $z$.  A direct
computation shows that the STFT of $e^{2\pi i \Phi_{2,z}(\zeta)}$ with
Gaussian window
 $\Psi(\zeta)=e^{-\pi|\zeta|^{2}}$ has modulus, up to constants,
 $e^{-\pi |u_{1}-u_{2}|^{2}/2}$, $u_{1},u_{2}\in\mathbb{R}^{2d}$.  This function belongs to $L^{\infty,1}_{v_{-4N}\otimes v_{s}}(\mathbb{R}^{{4d}})$ if and only if  $-4N+s<-2d$.
 This is in fact better than the condition $-2N+s<-2d$ in the assumptions,  due to the fact that here the derivatives of order $\geq 1$ of
  $\Phi_{2,z}(\zeta)$  are bounded from above by $\langle \zeta\rangle$, instead of by $\langle \zeta\rangle^{2}$, as in the general case.
  This explains the presence of $2N$ instead of $4N$. In the general case
  the best upper bound  for the derivatives of  $\Phi_{2,z}(\zeta)$ is
  $\langle \zeta\rangle^{2}$, so that the condition $-2N+s<-2d$
   should be sharp.  \par

\end{remark}
\begin{lemma}\label{modmult}
Let $m\in \mathcal{M}_v(\rd)$, $v$, $\nu, w$ be weight
 functions on $\rd$,  and $\{ f_z: z\in \rdd \}\subseteq \cS ' (\rd
 )$ be a set of distributions in $\cS '$.    If 
\begin{equation}\label{h1}
\sup_{z\in\rd}|V_\f f_z|\in L^{\infty,1}_{\nu^{-1} \otimes
m}(\rdd),
\end{equation}
for given  $\f\in\cS(\rd)$, then for every $h\in M^1_{w\nu\otimes
 v}(\rd)$  
\begin{equation}\label{h2}
\sup_{z\in\rd}|V_{\f ^2} (f_z h)|\in L^{1}_{w\otimes m}(\rdd) \, .
\end{equation}

\end{lemma}
\begin{proof} The proof of the multiplier property is similar to the
 convolution property of \modsp s in~\cite{CG}. 
Since
$V_gf\phas =(\widehat{f\cdot T_x \bar{g})}(\o)$ \cite[Lemma
3.1.1]{book}, 
we obtain the identity
$$
V_{\f ^2}  (f_z h ) (x,\eta ) = \big((f_z \overline{T_x\f} )( h T_x \overline{T_x\f})\big)\, \widehat{}
\, (\eta ) = (f_z \overline{T_x\f})\, \widehat{} \, \ast _\eta  (h
\overline{T_x\f}) \, \widehat{} \, (\eta)
\,  ,
$$
and consequently
$$
|V_{\f ^2} (f_z h)(x,\eta)|\leq |(\widehat{f_z\cdot T_x
\bar{\f})}|\ast_\o |(\widehat{h\cdot T_x \bar{\f})}|(\o) \, ,
$$
where the convolution is in the second variable $\eta $.
Now set $F(x,\eta ) = \sup _{z\in \rdd } |V_\f f_z (x,\eta )| = \sup
_{z\in \rdd } |(f_z\cdot T_x \bar{\f})\, \widehat{} \, (\eta ) | $ and
$H(x,\eta ) = |V_\f h (x,\eta )|$. Then
$$
\sup _{z\in \rdd } |V_{\f ^2} (f_z \, h)(x,\eta)| \leq (F \ast _\eta
H)(x,\eta ) \, .
$$
Finally,
\begin{align*}
 \|F \ast _\eta H \|_{L^1_{w\otimes m}} &= \intrd \intrd \intrd F(x,
   \eta - \zeta) H(x,\zeta ) d\zeta \, w(x) m(\eta ) \, dx d\eta \\
&\leq  \intrd \intrd \intrd F(x,
   \eta - \zeta) \nu (x)\inv m(\eta -\zeta )  H(x,\zeta )  \,
   w(x) \nu(x) v(\zeta ) \, dx d\eta d\zeta\\
&\leq \|F\|_{L^{\infty ,1}_{\nu \inv \otimes m} } \, \|H\|_{L^1_{\nu w
   \otimes v}} \, .
\end{align*}
In the last expression both norms are finite by assumption.
\end{proof}

\begin{corollary}\label{corollario} Let  $s\geq0$, $ N\in\mathbb{N},\quad N>\frac{s}{2}+d$ and $
  g\in M^1_{v_{{4N}}\otimes v_{{4N}}}(\rd)$ and $\Psi \in \cS (\rdd )$. 
Then the function $\Psi_z(w)=e^{2\pi i \Phi_{2,z}(w)}
\overline{g}\otimes\widehat{g}(w)$ defined as in \eqref{psi}
satisfies
$$\sup_{z\in\rdd}|V_{\Psi_z}\Psi (u_1,u_2)| \in L^1_{1\otimes v_{s}}(\R^{4d}).
$$
\end{corollary}
\begin{proof}
The symmetry property of the weight $\langle x \rangle ^{4N} \langle
\eta \rangle ^{4N}$ implies that  $M^1_{v_{{4N}}\otimes v_{{4N}}}(\rd )$
is invariant under the \ft , and thus $\hat{g}\in
M^1_{v_{4N}\otimes v_{4N}}(\rd)$. A tensor product argument then shows that
$\bar{g}\otimes\hat{g}\in M^1_{W}(\rdd)$ with
$W(x_{1},x_{2},\o_{1},\o_{2})=\langle x_{1}\rangle ^{{4N}}\langle
x_{2}\rangle ^{{4N}}\langle \o_{1}\rangle ^{{4N}}\langle \o_{2}\rangle
^{{4N}}$.  Since $W(x_1,x_2, \eta _1, \eta _2) \geq \langle (x_1,x_2)
\rangle ^{4N} \langle (\eta _1, \eta _2)\rangle ^{4N}$, we also obtain
that $\bar{g}\otimes\hat{g}\in M^1_{v_{4N}\otimes v_{4N}}(\rdd )$.

By Lemma~\ref{window} we have $\sup_{z\in\rdd}|V_\Psi
e^{2\pi i \Phi_{2,z}}|\in L^{\infty,1}_{v_{-4N}\otimes
v_{s}}(\R^{4d})$.

We now apply Lemma~\ref{modmult} with $f_z=  e^{2\pi i \Phi_{2,z}}$
and $h= \bar{g}\otimes\hat{g}$, $\nu = v_{4N}, m=v_s$, $w \equiv 1$,  and $v=v_{4N}$
(observe that $v_s$ is $v_{4N}$-moderate, since $s<4N$). As a
conclusion we obtain that $\sup_{z\in\rdd}|V_{\Psi}\Psi_z| \subseteq
L^1_{1\otimes v_{s}}(\R^{4d})$.
\end{proof}

We recall the following pointwise inequality of the \stft\
\cite[Lemma 11.3.3]{book}. It is often
useful when one needs to change window functions.
\begin{lemma}\label{changewind}
If  $g_0,g_1,\gamma\in\cS(\rd)$ such
that $\la \gamma, g_1\ra\not=0$ and
$f\in\cS'(\rd)$,  then the inequality
$$|V_{g_0} f\phas|\leq\frac1{|\la\gamma,g_1\ra|}(|V_{g_1} f|\ast|V_{g_0}\gamma|)\phas,
$$
holds pointwise for all $\phas\in\rdd$.
\end{lemma}
\begin{lemma}\label{wienerprop}
Let  $s\geq0$, $ N\in\mathbb{N},\quad N>\frac{s}{2}+d$ and $
  g\in M^1_{v_{{4N}}\otimes v_{{4N}}}(\rd)$.   Then
\begin{equation}\label{l1W}
\sup_{z\in\rdd}|V_{\Psi_z}\Psi| \in W(L^1_{1\otimes v_{s}})(\R^{4d}).
\end{equation}
\end{lemma}
\begin{proof} Let $\f\in\cS(\rdd)$ such that $\|\f\|_2=1$. Using Lemma \ref{changewind},
\begin{align*}
|V_{\Psi_z}\Psi|(u_1,u_2)&\leq|V_{\f}\Psi|\ast|V_{\Psi_z}\f|(u_1,u_2)\\
&\leq |V_{\f}\Psi|\ast (\sup_{z\in\rdd}|V_{\Psi_z}\f|)(u_1,u_2).
\end{align*}
Since $V_{\f}\Psi\in\cS(\R^{4d})\subset W(L^1_{1\otimes v_{s}})$ and
$\|\sup_{z\in\rdd}|V_{\Psi_z}\f | \, \|_{L^1_{1\otimes v_{s}}}<\infty$
by Corollary~\ref{corollario},  
the  convolution relation $L^1_{1\otimes v_{s}} \ast W(L^1_{1\otimes v_{s}}) \subseteq W(L^1_{1\otimes v_{s}})$ of \cite[Theorem
11.1.5]{book} implies that
$$\|\sup_{z\in\rdd}|V_{\Psi_z}\Psi|\,\|_{W(L^1_{1\otimes v_{s}})}\leq
\|V_{\f}\Psi\|_{W(L^1_{1\otimes v_{s}})}\|\sup_{z\in\rdd}|V_{\Psi_z}\f|\,\|_{L^1_{1\otimes v_{s}}}<\infty
\, ,
$$
and so  $\sup_{z\in\rdd}|V_{\Psi_z}\Psi |  \in W(L^1_{1\otimes v_{s}})$.
\end{proof}



We next formulate  an new  characterization of  the symbol classes  $M^{\infty,1}_{1\otimes
m}$ and $M^{\infty}_{1\otimes m}$ that is perfectly adapted to the
investigation of Fourier integral operators with a tame phase $\Phi $.

\begin{proposition}\label{proposition4}
Let $s\geq0$, $m\in \mathcal{M}_{v_{s}}(\rdd)$, $N$, $g$ be as in \eqref{g1} and $\sigma\in\cS'(\rdd)$. \\
(i) Then the symbol  $\sigma $ is in  $ M^{\infty,1}_{1\otimes
 m}(\R^{2d})$,  if and only if
\begin{equation}\label{p1}
\|\sup_{z\in\rdd}|V_{\Psi_z} \sigma|\|_{L^{\infty,1}_{1\otimes
    m}(\R^{4d})} = \int_{\rdd} \sup_{u_1\in\rd}\sup_{z\in\rdd}|V_{\Psi_z}
\sigma(u_1,u_2)|m(u_2)\,du_2<\infty,
\end{equation}
with $\|\sigma\|_{ M^{\infty,1}_{1\otimes m}}\asymp
\|\sup_{z\in\rdd}|V_{\Psi_z} \sigma|\|_{L^{\infty,1}_{1\otimes
m}}$.

In this case the function $H(u_2) =  \sup_{u_1\in\rd}\sup_{z\in\rdd}|V_{\Psi_z}
\sigma(u_1,u_2)|$ is in $W(L^1_m)(\rdd )$.

(ii) Likewise,  $\sigma \in M^{\infty}_{1\otimes m}(\R^{2d})$ if and only if
$\sup_{z\in\rdd}|V_{\Psi_z} \sigma|\in L^{\infty}_{1\otimes
m}(\R^{4d})$
with $\|\sigma\|_{ M^{\infty}_{1\otimes m}}\asymp
\|\sup_{z\in\rdd}|V_{\Psi_z} \sigma|\|_{L^{\infty}_{1\otimes
m}}$.
\end{proposition}
\begin{proof} We detail the proof of case (i). Case (ii) is
obtained similarly.  Let $\Psi\in\cS(\rdd)$  with
$\|\Psi\|_2=1$.\par

Assume first that  $\sigma \in M^{\infty,1}_{1\otimes
m}(\R^{2d})$. Then by Lemma \ref{changewind} we have
\begin{equation*}
|V_{\Psi_z}\sigma(u_1,u_2)|\leq|V_{\Psi}\sigma|\ast|V_{\Psi_z}\Psi|(u_1,u_2)\leq
|V_{\Psi}\sigma|\ast \sup_{z\in\rdd}|V_{\Psi_z}\Psi|(u_1,u_2).
\end{equation*}
Set $F(u_1,u_2)=\sup_{z\in\rdd}|V_{\Psi_z}\Psi(u_1,u_2)|$, so that
\[
|V_{\Psi_z}\sigma(u_1,u_2)|\leq \big(|V_{\Psi}\sigma|\ast F \big)(u_1,u_2).
\]
Finally
\begin{align*}
\|\sup_{z\in\rdd}|V_{\Psi_z} \sigma|\, \|_{L^{\infty,1}_{1\otimes m} }
& = \int_{\rdd}\sup_{u_1\in\rd}\sup_{z\in\rdd}|V_{\Psi_z}
\sigma(u_1,u_2)|m(u_2)\,du_2 \\
&\leq\| \, |V_{\Psi}\sigma | \ast
F\|_{L^{\infty,1}_{1\otimes m}}\\
&\leq \|V_{\Psi}\sigma\|_{L^{\infty,1}_{1\otimes
m}}\|F\|_{L^1_{1\otimes v_{s}}}\lesssim
\|\sigma\|_{M^{\infty,1}_{1\otimes m}}\|F\|_{L^1_{1\otimes
v_{s}}}<\infty,
\end{align*}
where in the last step we used the independence of the weighted $M^{\infty
 ,1}$-norm of the window~\cite[Thm.~11.3.7]{book} and Corollary \ref{corollario}. 

To obtain the sharper estimate, set $H(u_2) =  \sup_{u_1\in\rd}\sup_{z\in\rdd}|V_{\Psi_z}
\sigma(u_1,u_2)|$, $F_1(u_2) = \intrdd F(u_1,u_2) \, du_1$, and
$G_1(u_2) = \sup _{u_1\in \rdd } |V_\Psi \sigma (u_1, u_2)|$.
Then  the definition of $M^{\infty ,1}_{1\otimes m}(\rdd )$ implies
that $G_1 \in L^1_m(\rdd )$,  and Lemma~\ref{wienerprop}  implies that
$F_1\in W(L^1_{v_s})(\rdd )$. With these definitions we obtain a
pointwise estimate for $H$, namely
\begin{align*}
  H(w) & =   \sup_{u_1\in\rd}\sup_{z\in\rdd}|V_{\Psi_z}
\sigma(u_1,w)| \\
&\leq \intrdd \sup _{u_1} |V_\Psi \sigma (u_1,u_2)| F(z-u_1, w-u_2) du_1 du_2
\\
&= (G_1 \ast F_1)(w) \, .
\end{align*}
The convolution estimate \eqref{eq:may2} now yields that $H\leq G_1
\ast F_1 \in L^1_m \ast W(L^1_{v_s}) \subseteq W(L^1_{v_s})(\rdd )$,
as claimed.

\par Conversely, assume
\eqref{p1}. Using  Lemma \ref{changewind} again, we deduce that
\[
|V_{\Psi}\sigma (u_1,u_2)|\leq\frac{1}{\langle \Psi_z,\Psi_z\rangle}|V_{\Psi_z}\sigma|\ast|V_{\Psi}\Psi_z|(u_1,u_2),
\]
and
\begin{align*}
\langle \Psi_z,\Psi_z\rangle&=\int_\rdd |e^{2\pi
 i\Phi_{2,z}(w)}(\overline{g}\otimes\widehat{g})(w)|^2\,dw\\
&=\int_{\rdd} |(\overline{g}\otimes\widehat{g})(w)|^2\,dw = \|\overline{g}\|^2_2\|\widehat{g}\|^2_2=\|g\|_2^4
\end{align*}
is in fact a constant (depending on $g$). Using the involution $
f^*(z) = \overline{f(-z)}$, we continue with
\begin{align*}
|V_{\Psi}\sigma|(u_1,u_2)&\lesssim  |V_{\Psi_z}\sigma|\ast|V_{\Psi}\Psi_z|(u_1,u_2)\\
&=|V_{\Psi_z}\sigma|\ast|(V_{\Psi_z}\Psi)^\ast|(u_1,u_2)\\
&\leq \left(\sup_{z\in\rdd} |V_{\Psi_z}\sigma |\right)\ast \left( \sup_{z\in\rdd} |(V_{\Psi_z}\Psi)^\ast| \right)(u_1,u_2).
\end{align*}
Since by Corollary~\ref{corollario}
\[
\sup_{z\in\rdd} |(V_{\Psi_z}\Psi)^\ast|=\sup_{z\in\rdd} |V_{\Psi_z}\Psi|\in L^{1}_{1\otimes v_{s}}(\R^{4d}),
\]
we get
\[
\|V_{\Psi}\sigma\|_{L^{\infty,1}_{1\otimes m}}\lesssim
\|\sup_{z\in\rdd}|V_{\Psi_z}\sigma|\|_{L^{\infty,1}_{1\otimes m}}\|
\sup_{z\in\rdd} |V_{\Psi_z}\Psi|\|_{L^{1}_{1\otimes v_{s}}}<\infty \, ,
\]
and the proof is complete. 
\end{proof}

We  now  prove Theorems \ref{6bis} and ~\ref{6biss}.
\begin{proof}[Proof of Theorem \ref{6bis}]
We know from \cite[Equ.~(39)]{fio1} or  \cite[Equ.~(3.3)]{fio5} that
\begin{align}
|\langle T\pi(x,\o)g,\pi(x',\o')g
\rangle|&=|V_{\Psi_{(x',\o)}}\sigma
\big((x',\o),(\o'-\nabla_x\Phi(x',\o),x-\nabla_\eta\Phi(x',\o))\big)|\notag\\
&\leq \sup_{z\in\rdd} |V_{\Psi_{z}}\sigma
\big((x',\o),(\o'-\nabla_x\Phi(x',\o),x-\nabla_\eta\Phi(x',\o))\big)|.\label{stima}
\end{align}

We claim that the function that controls the off-diagonal decay of
\linebreak  $|\langle T\pi(x,\o)g,\pi(x',\o')g
\rangle|$ is exactly the function $ H(w):= \sup _{u_1 \in \rdd } \sup_{z\in\rdd}
|V_{\Psi_{z}}\sigma (u_1,w)|$ introduced in
Proposition~\ref{proposition4}. There   we have already proved that
$H$ is in $W(L^1_{m})(\rdd)$,  and now  \eqref{stima} implies that
$$
|\langle T \pi(x,\o)g,\pi(x',\o')g\rangle|\leq
H(\o'-\nabla_x\Phi(x',\o),x-\nabla_\eta\Phi(x',\o)).
$$
\end{proof}


A decay
estimate in terms of the canonical transformation $\chi$ can be
obtained by imposing a stronger condition on the symbol class. 

\smallskip
\begin{proof}[Proof of Theorem \ref{6biss}]
We follow the pattern of the proof of Theorem \ref{6bis}.
If $\Phi $ is a tame phase function, then the argument of $V_{\Psi
  _z}\sigma$ in \eqref{stima} can be estimated further as
\begin{equation}
  \label{eq:chh8}
|\o'-\nabla_x\Phi(x',\o) + |x-\nabla_\eta\Phi(x',\o)| \geq
|\o'-\chi_2\phas | + | x'-\chi_1\phas |
\end{equation}
by \cite[Lemma 3.1]{fio1}.

Set $H(u_1,u_2):=\sup_{z\in\rdd} |V_{\Psi_{z}}\sigma (u_1,u_2)|$.
By Proposition \ref{proposition4} (ii),
$H\in L^{\infty}_{1\otimes v_s}(\R^{{4d}})$ and
$\|H\|_{L^{\infty}_{ 1\otimes v_s}}\lesssim \|\sigma\|_{
M^{\infty}_{1\otimes v_s}}$. Now, using \eqref{stima},
\begin{align*}
|\langle  T\pi(x,\o)g,\pi(x',\o')g \rangle|&\leq \sup_{u_1\in\rdd}
H(u_1,\o'-\nabla_x\Phi(x',\o),x-\nabla_\eta\Phi(x',\o))\\
&\lesssim \sup_{u_1\in\rdd}
H(u_1,\o'-\nabla_x\Phi(x',\o),x-\nabla_\eta\Phi(x',\o))\\
&\quad\quad\quad\times\,\frac{\la\o'-\nabla_x\Phi(x',\o),x-\nabla_\eta\Phi(x',\o)\ra^s}{\la
\o'-\nabla_x\Phi(x',\o),x-\nabla_\eta\Phi(x',\o)\ra^s}\\
&\lesssim
\frac{\sup_{u_2} \sup_{u_1\in\rdd} H(u_1,u_2)v_s(u_2)}{\la
\o'-\nabla_x\Phi(x',\o),x-\nabla_\eta\Phi(x',\o)\ra^s}\\
&= \frac{\|H\|_{L^\infty_{1\otimes v_s}}}{\la
\o'-\nabla_x\Phi(x',\o),x-\nabla_\eta\Phi(x',\o)\ra^s}\\
&\lesssim\frac{\|\sigma\|_{M^{\infty}_{1\otimes v_s}}}{\la
\o'-\chi_2\phas,x'-\chi_1\phas\ra^s},
\end{align*}
where in the last inequality we have  used  \eqref{eq:chh8}.
\end{proof}

\section{Gabor multipliers}\label{gaborsection}
We now introduce a type of Gabor multipliers that are  tailored to the 
canonical transformation $\chi$ of a FIO. In the
case of the identity map $\chi={\rm id}$ we get the standard Gabor
multipliers studied,  e.g.,  in \cite{Fei-Now02,approx}. \par Since
in general the map $\chi$  does not preserve a given lattice
$\Lambda$, we first replace $\chi $ by a mapping $\chi'$ that enjoys
$\chi'(\Lambda)\subset\Lambda$.\par

We define the integer part of a vector
$y=(y_1,...,y_{2d})$ as
\[
\lfloor y\rfloor=(\lfloor
y_1\rfloor,...,\lfloor
y_{2d}\rfloor),
\]
i.e., by taking the integer
parts of its components.
\par
Now consider a lattice of the form $\Lambda=A\zdd$ with $A\in GL(2d,\R)$.
Given a  map $\chi:\R^{2d}\to\R^{2d}$, we approximate $\chi$ by a
map $\chi':\Lambda\to\Lambda $
defined by
\[
\chi'(\lambda)=A\lfloor
A^{-1}\chi(\lambda)\rfloor.
\]
Observe that $|y-\lfloor
y\rfloor|<\sqrt{2d}$ and therefore
\[
|\chi(\lambda)-\chi'(\lambda)|<\sqrt{2d}\|A\|\, .
\]
The almost diagonalization of  Theorem
\ref{6biss} can now be formulated in terms of $\chi '$ as follows:
using the
inequality $\l x\r\leq
\sqrt{2}\l x+y\r \l y \r$ we
have
\begin{equation}\label{5fbis}
   |\la T \pi(\mu)g, \pi(\lambda)g\ra|\leq C_N 2^s
   \l \sqrt{2d}\,\|A\|\r^{2s}\,{\la \chi'(\mu)
   -\lambda\ra^{-2s}
   },\qquad \forall \lambda,\mu\in\Lambda.
\end{equation}
We collect here some properties which will be used in the sequel.
\begin{lemma}\label{proposition4.1} Let $v(z)=v_s(z)=\langle z\rangle^s$, $z\in\rdd$, $s\geq0$, and $m\in\mathcal{M}_v(\rdd)$. Then
\begin{equation}\label{equa6}
v\circ\chi\asymp v,\qquad v\circ\chi'\asymp v,
\end{equation}
and
\begin{equation}\label{equa7}
m\circ \chi\in\mathcal{M}_v(\rdd)\,  , \qquad m\circ \chi'\in
\mathcal{M}_v(\rdd) .
\end{equation}
\begin{proof}
Since $\chi$ is Lipschitz continuous,  we have
\[
|\chi(z)|\leq |\chi(0)|+|\chi(z)-\chi(0)|\leq|\chi(0)|+C|z|\leq
C'(1+|z|) \qquad \forall z\in \rdd \, .
\]
By applying the same argument to $\chi^{-1}$, which is also  Lipschitz
continuous, we obtain  $1+|\chi(z)|\asymp 1+|z|$. Since we also
have $\chi'(z)=\chi(z)+\mathcal{O}(1)$, both the estimates in
\eqref{equa6} follow.\par
Concerning \eqref{equa7}, observe that
\[
m(\chi(z+\zeta))=m(\chi(z+\zeta)-\chi(z)+\chi(z))\leq v(\chi(z+\zeta)-\chi(z))m(\chi(z))\lesssim v(\zeta)m(\chi(z)),
\]
where in the last step we used the fact that
$|\chi(z+\zeta)-\chi(z)|\lesssim |\zeta|$. This proves the
first formula in \eqref{equa7}. The second formula is obtained
similarly, because $\chi'(z)=\chi(z)+\mathcal{O}(1)$, so that
$\chi'(z+\zeta)-\chi'(z)=\chi(z+\zeta)-\chi(z)+\mathcal{O}(1)$,
and we still have  $1+|\chi'(z+\zeta)-\chi'(z)|\lesssim
1+|\zeta|$.
\end{proof}
\begin{remark}\rm
Notice that the above proposition holds only for weights with  polynomial
growth,  but not for weights with super-polynomial growth, e.g.,
$v(z)=e^{a|z|^b}$, $0<b<1$, $a>0$.  This is related to
the fact that the weights $v_s$ satisfy the doubling condition $v_s
(2z) \asymp v_s(z)$, whereas  weights with super-exponential growth do
not. 
\end{remark}
\end{lemma}
Now, let $\G(g,\Lambda)$ be a Gabor system, and $\chi$, $\chi'$ be the
canonical transformations of an FIO. Given a sequence ${\bf a}=\left({
a}_\nu\right)_{\nu\in\Lambda}$, we define (formally) the Gabor
multiplier
\[
M_{\bf a} f=M_{\bf a}^{\chi',g,\Lambda}f=\sum_{\lambda\in\Lambda}
a_\lambda \, \langle f,\pi(\lambda) g\rangle\pi(\chi'(\lambda))g.
\]
In order to give a precise meaning to this definition, we need to
study the convergence of the above series.

\begin{lemma}\label{lemma2.2}
Let 
$g\in M^1_v(\rd)$, $m\in\mathcal{M}_v(\rdd)$,  and $1\leq
p\leq\infty$. If $h = (h_\lambda)_{\lambda \in \Lambda } \in \ell
^p_{m\circ \chi '} (\Lambda )$, then
\begin{equation}
 \label{eq:new2}
\|  \sum_{\lambda\in\Lambda}h_\lambda \pi(\chi'(\lambda))
g\|_{M^p_{m}}\leq C\Big(\sum_{\lambda\in \Lambda}|h_\lambda|^p
m(\chi'(\lambda))^p\Big)^{1/p} \, .
\end{equation}
The series $\sum_{\lambda\in\Lambda}h_\lambda \pi(\chi'(\lambda))
g $ converges unconditionally in $M^p_m(\rd)$ for $p<\infty $ and
weak$^*$ unconditionally for $p=\infty $.
\end{lemma}
\begin{proof}
We assume  $p<\infty$ and leave the  modification for  the case
$p=\infty$ to the reader.

Since the finite sequences are norm-dense in $\ell ^p_m(\Lambda)$
for $p< \infty$ (and weak$^*$-dense for $p=\infty$), it suffices
to show \eqref{eq:new2} for finite sequences. For a  finite set
$F\subset \Lambda$ we set
$\tilde{F}=\chi'(F)$.  Since $g\in M^1_v(\rd)$, 
the synthesis operator $\{c_\alpha\}\mapsto\sum_{\alpha\in\Lambda}
c_\alpha\pi(\alpha) g$ is bounded from $\ell ^p_{m}(\Lambda)$ to
$M^p_m(\rd)$~\cite[Thm.~12.2.4]{book}, so that
\begin{align}
\|\sum_{\lambda\in F}h_\lambda \pi(\chi'(\lambda)) g\|_{M^p_m}
&=\|\sum_{\alpha\in\tilde{F}}\Big(\sum_{\lambda\in
 {\chi'}^{-1}(\alpha)}h_\lambda\Big)\pi(\alpha) g\|_{M^p_m} \notag\\
&\lesssim \Big(\sum_{\alpha\in\tilde{F}}\Big|\sum_{\lambda\in
 {\chi'}^{-1}(\alpha)}h_\lambda\Big|^p m(\alpha)^p\Big)^{1/p}. \label{eq:new1}
\end{align}
Since $\chi \inv $ is  Lipschitz
continuous, $\chi'$ is ``uniformly almost injective'', in the sense that
$\sup_{\alpha\in\Lambda}\#
{\chi'}^{-1}(\{\alpha\})<\infty$. Hence the last expression in
\eqref{eq:new1}  is bounded
by
\[
\lesssim\Big(\sum_{\alpha\in\tilde{F}}\sum_{\lambda\in {\chi'}^{-1}(\alpha)}|h_\lambda|^p m(\alpha)^p\Big)^{1/p}\\
=\Big(\sum_{\lambda\in F}|h_\lambda|^p m(\chi'(\lambda))^p\Big)^{1/p}.
\]

By density \eqref{eq:new2} holds for all $(h_\lambda ) \in \ell
^p_{m\circ \chi '}(\Lambda)$. The unconditional convergence
follows from the
norm estimate. 
\end{proof}
\begin{proposition}\label{prop4}
Let 
$g\in M_{v^2} ^1(\rd)$,   $m,\tilde{m}\in\mathcal{M}_v(\rdd)$, and
$1\leq p\leq\infty$. If ${\bf a}=(a_\lambda)_{\lambda\in\Lambda}$
is a
sequence in $\ell ^\infty _{\tilde{m}}(\Lambda)$, 
then the Gabor multiplier $M_{\bf a}=M_{\bf a}^{\chi',g,\Lambda}$
is bounded from $M^{p}_{\frac{m\circ\chi'}{\tilde{m}}}(\rd)$ to
$M^p_m(\rd)$. Its operator norm is bounded by  $\|M_{\bf
  a}\|_{M^{p}_{\frac{m\circ\chi'}{\tilde{m}}}\to M^p_m}\lesssim \|{\bf a}\|_{\ell
  ^\infty_{\tilde{m}}}$.
\end{proposition}
\begin{proof}
Since $\frac{m\circ\chi'}{\tilde{m}} \in \cM _{v^2}(\rd)$ and
$g\in
M^1_{v^2} (\rd)$, 
the coefficient operator  $f\mapsto\langle f,\pi(\lambda)g\rangle$
is bounded from $M^{p}_{\frac{m\circ\chi'}{\tilde{m}}}(\rd)$ to
$\ell ^{p}_{\frac{m\circ\chi'}{\tilde{m}}}(\Lambda)$,
see~\cite[Thm.~12.2.3.]{book}. Hence, if  ${\bf
 a}=(a_\lambda)_{\lambda\in\Lambda}\in \ell
^\infty_{\tilde{m}}(\Lambda)$ and  $f\in
M^{p}_{\frac{m\circ\chi'}{\tilde{m}}}(\rd)$, then  the sequence
$(a_\lambda \langle f,\pi(\lambda)g\rangle)_{\lambda\in\Lambda}$
belongs to $\ell ^{p}_{m\circ\chi'}(\Lambda)$. Therefore, by Lemma
\ref{lemma2.2} the series $\sum_{\lambda\in\Lambda}
a_\lambda\langle f,\pi(\lambda) g\rangle \pi(\chi'(\lambda))g$
converges unconditionally in $M^p_m(\rd)$. The estimate for the
operator norm follows from
\begin{align*}
\|\sum_{\lambda\in\Lambda} a_\lambda\langle f,\pi(\lambda) g\rangle \pi(\chi'(\lambda))g\|_{M^p_m}&\lesssim \Big(\sum_{\lambda\in \Lambda}|a_\lambda\langle f,\pi(\lambda) g\rangle|^p m(\chi'(\lambda))^p\Big)^{1/p}\\
&\leq\|\tilde{m}a\|_{\ell ^\infty}  \Big(\sum_{\lambda\in\Lambda}|\langle f,\pi(\lambda)g\rangle|^p\Big(\frac{m(\chi'(\lambda))}{\tilde{m}(\lambda)}\Big)^p\Big)^{1/p}\\
&\lesssim \|a\|_{\ell ^\infty_{\tilde{m}}}\|f\|_{M^p_{\frac{m\circ\chi'}{\tilde{m}}}}.
\end{align*}
\end{proof}

\begin{cor} \label{osservazione}
 If ${\bf a } \in \ell ^\infty _{v_s}(\Lambda)$ for $s\in \bR$ and $g\in \cS
 (\rd )$, 
then $M_{\bf a}^{\chi',g,\Lambda }$ maps $\cS(\rd)$ into $\cS(\rd)$.
\end{cor}
\begin{proof}
Since   $v_N\circ\chi'\asymp v_N$ by Lemma~\ref{proposition4.1},
$M^p_{v_N \circ \chi' / v_s} (\rd) = M^p_{v_{N-s}}(\rd)$. By 
Proposition \ref{prop4} $M_{\mathbf{a}}$ maps $M^p_{v_{N-s}}(\rd)$
to $M^p_{N}(\rd)$. Since $\cS(\rd)=\bigcap_{N\geq0} M^\infty
_{v_N}(\rd)$ by \eqref{osservazione0}, $M_{\mathbf{a}}$ maps $\cS
(\rd ) = \bigcap _{N\geq 0}  M^\infty_{v_{N-s}}(\rd)$ into $ \cS (\rd
) $.
\end{proof}

\section{Approximation of FIOs by Gabor Multipliers}
Our next aim is to approximate a FIO $T$ with a tame phase  by Gabor
multipliers associated to a Parseval frame
$\G(g,\Lambda)$. First we will derive a representation
\begin{equation}
  \label{jaha}
T=\sum_{\nu\in\Lambda}\pi(\nu) M_{{\bf a}_\nu},  
\end{equation}
with convergence in several operator norms.  To find  the candidate
symbols ${\bf a}_\nu$, we  argue as in~\cite{approx}.
\par Using the commutation relations $M_\eta T_x=e^{2\pi i
x \eta}T_xM_\eta$, we can write
\begin{equation}\label{symb0}
\pi(\chi'(\mu)+\nu)=c_{\nu,\mu}\pi(\nu)\pi(\chi'(\mu)), \quad
\text{ for } \nu,\mu\in\Lambda \, ,
\end{equation}
with $|c_{\nu,\mu}|=1$. We will  show that the choice of
\begin{equation}\label{symb}
{\bf a}_\nu(\mu)=c_{\nu,\mu}\la T\pi(\mu) g,\pi(\chi'(\mu)+\nu)g\ra
\end{equation}
leads to the formal  representation~\eqref{jaha}.

Since $\mathcal{G}(g,\Lambda)$ is a Parseval frame, $f$ and $Tf$
possess the  Gabor expansions
$$f=\sum_\mu\la f,\pi(\mu)g\ra\pi(\mu)g,\quad \text{ and } \quad
Tf=\sum_\lambda\la Tf,\pi(\lambda)g\ra\pi(\lambda) g.
$$
Setting $\nu=\lambda-\chi'(\mu)\in \Lambda $,  we can write
\begin{align}
               Tf&=\sum_\mu\la f, \pi(\mu)g\ra\sum_\lambda \Big\la
               T\pi(\mu)g,\pi(\lambda)g\Big\ra\pi(\lambda)g \notag \\
               &=\sum_\mu\sum_\nu\la f, \pi(\mu)g\ra\, \Big\la
               T\pi(\mu)g,\pi(\chi'(\mu)+\nu)g\Big\ra\pi(\chi'(\mu)+\nu)
               g \notag\\
               &=\sum_\mu\sum_\nu c_{\nu,\mu}\la f, \pi(\mu)g\ra \,
               \Big\la T\pi(\mu)g,\pi(\chi'(\mu)+\nu)g \Big\ra  \,
               \pi(\nu)\pi(\chi'(\mu))g \notag \\
               &=\sum_\nu\pi(\nu)\sum_\mu c_{\nu,\mu}\la f,
               \pi(\mu)g\ra\la
               T\pi(\mu)g,\pi(\chi'(\mu)+\nu)g\ra\pi(\chi'(\mu))g
               \notag \\
               &=\sum_\nu\pi(\nu)M_{{\bf a}_\nu}f. \label{eq:haha26}
\end{align}
The following result  gives a precise meaning to the above
computation for test functions in the  Schwartz class $\cS (\rd ) $.
\begin{proposition}\label{proposition5.1}
Let $\G(g,\Lambda)$ be a Parseval frame and
$g\in\cS(\rd)\setminus\{0\}$.\par\medskip\noindent (i) If $|{\bf
a}_\nu(\mu)|\leq C(1+|\mu|+|\nu|)^N$
for some $N\geq0$, then the series $\sum_{\nu\in\Lambda}\pi(\nu)
M_{a_\nu}$ converges unconditionally  in the strong  topology of
$\mathcal{L}(\cS(\rd),\cS'(\rd))$ and defines a continuous
operator from $\cS(\rd)$ to $\cS'(\rd)$ ($\cS'(\rd)$ is always
endowed with the weak$^\ast$ topology).\par\medskip\noindent (ii)
Let $A$ be a continuous operator from $\cS(\rd)$ to $\cS'(\rd)$,
and choose ${\bf a}_\nu(\mu)$, $\nu,\mu\in\Lambda$ as
\[
{\bf a}_\nu(\mu)=c_{\nu,\mu}\la A\pi(\mu) g,\pi(\chi'(\mu)+\nu)g\ra
\]
with the constant $c_{\nu,\mu}$ as in \eqref{symb0}. Then the grow estimates
\begin{equation}\label{grow1}
|{\bf a}_\nu(\mu)|\leq C(1+|\mu|+|\nu|)^N,
\end{equation}
are satisfied  for some  constants $C,N\geq 0$  depending on $A$ and
$g$.\par Furthermore, $A$ can be represented as the sum
$\sum_{\nu\in\Lambda}\pi(\nu) M_{a_\nu}$ of shifted Gabor
multipliers in the strong operator topology of
$\mathcal{L}(\cS(\rd),\cS'(\rd))$.
\end{proposition}
\begin{proof}
(i) By Corollary~\ref{osservazione} the Gabor multiplier
\[
M_{\mathbf{a}_\nu}f=\sum_{\mu\in\Lambda}{\bf a}_\nu(\mu)\langle
f,\pi(\mu)g\rangle \pi(\chi'(\mu))g\in\cS(\rd)
\]
maps $\cS (\rd ) $ to $\cS (\rd ) $. We need to show that  the series
$\sum_{\nu\in\Lambda} \pi(\nu) M_{{\bf a}_\nu}f$, $f\in \cS (\rd )$,
converges unconditionally in $\cS '(\rd)$. 
\par Let  $h\in\cS(\rd)$ and choose $s$ and $l$ large enough, precisely
 $s>N+2d$ and  $l> N+s+2d$. Then 
\begin{align}\label{equa1}
|\langle \sum_\nu\pi(\nu)M_{{\bf a}_\nu}f,h\rangle|&=|\sum_\nu\pi(\nu)\sum_{\mu} {\bf a}_\nu(\mu)\langle f,\pi(\mu)g\rangle \pi(\chi'(\mu))g,h\rangle|\nonumber\\
&\leq \sum_\nu\sum_\mu|{\bf a}_\nu(\mu)||\langle f,\pi(\mu)g\rangle||\langle \pi(\nu+\chi'(\mu))g,h\rangle|\nonumber\\
&\lesssim  \sum_\nu\sum_\mu (1+|\mu  | + |\nu |)^N \,
(1+|\mu|))^{-l}(1+|\chi'(\mu)|)^s(1+|\nu|)^{-s}<\infty.
\end{align}
In the last step we used the assumption $|a_\nu(\mu)|\leq
C(1+|\mu|+|\nu|)^N$ and the fact that $1+|\chi'(\mu)|\asymp
1+|\mu|$ from Lemma~\ref{proposition4.1}. Then  the double series
converges because of our choice   $s>N+2d$ and  $l> N+s+2d$.
\par\medskip (ii) The proof is similar
to that of Proposition 5 (ii) of \cite{approx}. Let
$K\in\cS'(\rdd)$ be the Schwartz kernel of the operator $A$.
Writing $\mu=(\mu_1,\mu_2)$, $\nu=(\nu_1,\nu_2)\in \bR ^{4d}$ and
$\chi'(\mu)=(\chi'_1(\mu),\chi'_2(\mu))\in \bR ^{4d} $,  we have
\begin{align*}
|{\bf a}_\nu(\mu)|&=|\langle A\pi(\mu) g,\pi(\chi'(\mu)+\nu) g\rangle|\\
&=|\langle K, \pi(\chi'(\mu)+\nu) g\otimes \overline{\pi(\mu)g}\rangle|\\
&=| V_{g\otimes \overline{g}} K(\chi'_1(\mu)+\nu_1,\mu_1,\chi'_2(\mu)+\nu_2,-\mu_2)|\\
&\leq C(1+|\chi'(\mu)|+|\mu|+|\nu|)^N
\end{align*}
for some $C,N>0$.  In the last step we used the fact that the STFT of
a tempered distribution grows at most polynomially. Since
$1+|\chi'(\mu)|\lesssim 1+|\mu|$ we see that the sequence ${\bf
  a}_\nu(\mu)$ satisfies the polynomial growth
condition~\eqref{grow1}.  It follows then from part  (i) that the
series $\sum_{\nu\in\Lambda}\pi(\nu) M_{{\bf a}_\nu}$ converges in
$\mathcal{L}(\cS(\rd),\cS'(\rd))$. Furthermore, its sum must
coincide with $A$, because the formal   computation
\eqref{eq:haha26} was justified  in part (i) and therefore
\[
\langle \sum_\nu\pi(\nu)M_{{\bf a}_\nu}f,h\rangle=\langle A f,h\rangle,\qquad \forall f,h\in\cS(\rd).
\]
(One needs to interchange two summations, which is possible by \eqref{equa1}).
\end{proof}
\begin{remark}\label{osservazione5}\rm
Proposition~\ref{proposition5.1} extends to more general windows
provided that  the sequences ${\bf a_\nu}$ satisfy stronger
estimates. \par For example, if $\sup _{\mu , \nu \in \Lambda
}|{\bf a_\nu}(\mu)|= C<\infty $,  then one may use  windows $g\in
M^1(\rd)$. In fact, for
$f,h\in M^1(\rd)$, \eqref{equa1} can be modified to yield
\begin{align*}
|\langle \sum_\nu\pi(\nu)M_{{\bf a}_\nu}f,h\rangle| 
&\leq \sum_\nu\sum_\mu|{\bf a}_\nu(\mu)||\langle
f,\pi(\mu)g\rangle||\langle \pi(\nu+\chi'(\mu))g,h\rangle|\\
&\leq C\sum_{\mu} |\langle f,\pi(\mu)g\rangle| \sum_\nu  |\langle
\pi(\nu+\chi'(\mu))g,h\rangle| \\
&\leq C' \|f \|_{M^1} \, \|h \|_{M^1} <\infty \, .
\end{align*}
where in the last step we have used the 
characterization \eqref{osservazione0} of $M^1(\rd )$ (with $m\equiv1$ and $p=1$).
\end{remark}

We now consider the convergence of series of Gabor multipliers on
modulation spaces. 
In what follows the space ${\mathcal M}^\infty_m(\rd)$ denotes the
closure of the Schwartz class with respect to the $
M^\infty_m$-norm. Whereas Gabor expansions converge only weak$^*$
on $M^\infty _m(\rd )$, they are norm convergent of $\mathcal{M}
^\infty _m(\rd )$.

\begin{proposition}\label{proposition3}
Let $\G(g,\Lambda)$ be a Parseval frame and $g\in M^1_v(\rd)$. If
the sequence of symbols ${\bf a}_\nu$ satisfies
\[
\sum_{\nu\in\Lambda}\|{\bf a}_\nu\|_{\ell ^\infty}v(\nu)<\infty \, ,
\]
then  the series $\sum_\nu\pi(\nu)M_{{\bf a}_\nu}$ converges in
$\mathcal{L}(M^p_{m\circ\chi'},M^p_m)$ for every $p\in[1,\infty )$ and
$v$-moderate weight $m$. If $p=\infty $, the series
$\sum_\nu\pi(\nu)M_{{\bf a}_\nu}$ converges in
$\mathcal{L}(\mathcal{M}^\infty
_{m\circ\chi'},\mathcal{M}^\infty_m)$.
\end{proposition}
\begin{proof}
By Proposition \ref{prop4} each operator $M_{{\bf a}_\nu}$ is
bounded from  $M^{p}_{m\circ\chi'}(\rd)$ to $ M^p_m(\rd)$ with the
norm being  dominated by $\|{\bf a}_\nu\|_{\ell ^\infty}$.
Furthermore, $\|\pi(\nu)\|_{M^p_m\to M^p_m}\lesssim v(\nu)$. Hence
\begin{align*}
\sum_\nu\| \pi(\nu)M_{{\bf a}_\nu}\|_{M^{p}_{m\circ\chi'}\to M^p_m}&
\leq \sum_\nu \|\pi(\nu)\|_{M^p_m\to M^p_m}\|M_{{\bf a}_\nu}\|_{M^p_{m\circ\chi'}\to M^p_m}\\
&\lesssim \sum_\nu\|{\bf a}_\nu\|_{\ell ^\infty}v(\nu)<\infty.
\end{align*}
This gives the desired conclusion.
\end{proof}

We return to the study of Fourier integral operators. The following
approximation theorem is the main result of our work. 

\begin{theorem}\label{teorema2.1}
Let $\G(g,\Lambda)$ be a Parseval frame and $g\in \cS(\rd)$.
Let $T $ be a Fourier integral operator with a tame phase and with the
associated lattice transformation $\chi ' $, and define the multiplier
symbol as
$$
{\bf a}_\nu(\mu)=c_{\nu,\mu}\la T\pi(\mu) g,\pi(\chi'(\mu)+\nu)g\ra \,
.
$$

(i)  If  $\sigma \in W^\infty _{2N}(\rdd )$ for $N>d$ and if $0\leq r<
2N-2d$, then the
series $\sum_{\nu\in\Lambda} \pi(\nu)M_{{\bf a}_\nu}$ converges to
$T$ in $\mathcal{L}(M^p_{m\circ\chi'},M^p_m)$ and in
$\mathcal{L}({\mathcal M}^\infty_{m\circ\chi'},{\mathcal
M}^\infty_m)$, for every $p\in[1,\infty)$,  and $v_r$-moderate
weight $m$ ($v_r(z)=\langle z\rangle^r$, $z\in\rdd$), and
\[
\|T-\sum_{|\nu|\leq L} \pi(\nu)M_{\mathbf{a}_\nu}\|_{M^p_{m\circ\chi'}\to
  M^p_m}\lesssim L^{r+2d-2N}.
\]
(ii) If the symbol $\sigma$ is in $M^\infty_{1\otimes v_s}(\rdd)$
with $s>2d$, then $T=\sum_{\nu \in \Lambda } \pi(\nu)M_{\mathbf{a}_\nu}$
with the error estimate
$$
\|T-\sum_{|\nu|\leq L} \pi(\nu)M_{\mathbf{a}_\nu}\|_{M^p_{m\circ\chi'}\to
  M^p_m}\lesssim L^{r+2d-s}\, .
$$
\end{theorem}
\begin{proof}
(i) 
We first estimate the magnitude of the multiplier  symbols
$\mathbf{a}_\nu $:  Theorem~\ref{old}, written with $\chi '$ as
in \eqref{5fbis},  implies that
\begin{align*}
|{\bf a}_\nu(\mu)|=|\la
T\pi(\mu)g,\pi(\chi'(\mu)+\nu)g\ra|\lesssim
\l
\chi'(\mu)-(\chi'(\mu)+\nu)\ra^{-2N}=\la\nu\ra^{-2N},
\end{align*}
for every $\mu\in\Lambda$. Hence
$$ \|\mathbf{a}_\nu \|_\infty \lesssim \langle \nu \rangle ^{-2N} \,
$$ and
\[
\sum_{\nu\in\Lambda}\langle \nu\rangle^r\|{\bf
a}_\nu\|_{\ell ^\infty}\leq \sum_{\nu\in\Lambda}\langle
\nu\rangle^{r-2N}<\infty \, .
\]
By Proposition \ref{proposition3}  the series $\sum_{\nu}
\pi(\nu)M_{\mathbf{a}_\nu}$ converges in the operator norm from
$M^p_{m\circ\chi'}(\rd)\to M^p_m(\rd)$  for every $1\leq p\leq
\infty$. Since the series converges to $T$ in
$\mathcal{L}(\cS(\rd),\cS'(\rd))$ by
Proposition \ref{proposition5.1}, the sum must be identical to  $T$. 
\par For the  error estimate we
observe that
\begin{align*}
\|T-\sum_{|\nu|\leq L}
\pi(\nu)M_{\mathbf{a}_\nu}\|_{M^p_{m\circ\chi'}\to M^p_m}&\leq
\sum_{|\nu|>L}\|\pi(\nu)\|_{M^p_m\to M^p_m}
\|M_{\mathbf{a}_\nu}\|_{M^p_{m\circ\chi'}\to M^p_m}\\ 
&\lesssim \sum_{|\nu|>L} \langle\nu\rangle^r\|\mathbf{a}_\nu\|_{\ell^\infty}\\
&\lesssim \sum_{|\nu|> L}\la\nu\ra^{r-2N}\lesssim L^{r+2d-2N}.
\end{align*}

(ii) is proved similarly by using the decay estimate of Theorem~\ref{6biss}
instead of Theorem~\ref{old}.
\end{proof}

As  a byproduct of Theorem~\ref{teorema2.1} we obtain
an alternative proof of  \cite[Theorem 4.1]{fio1}:
\begin{corollary} Under the assumptions of Theorem \ref{teorema2.1},
 the  Fourier integral operator $T$ is a bounded
operator from $M^p_{m\circ\chi'}(\rd)$ to $M^p_m(\rd)$, simultaneously
for every
$1\leq p< \infty$,  and from ${\mathcal
M}^\infty_{m\circ\chi'}(\rd)$ to ${\mathcal M}^\infty_m(\rd)$.
\end{corollary}

\begin{remark}
  So far we have used without loss of generality that $\cG (g,
  \Lambda ) $ is a tight Gabor frame.  If $\cG (g,\Lambda )$ is an
  arbitrary frame with $g\in M^1_v(\rd)$, then there exists a dual window
  $\gamma \in M^1_v(\rd)$ such that every $f\in M^p_m(\rd)$, for $1\leq p\leq
  \infty $ and $m\in \cM _v(\rdd)$,  possesses the Gabor expansion $f= \sum
  _{\lambda \in \Lambda } \langle f, \pi (\lambda ) g\rangle \pi
  (\lambda ) \gamma $ with convergence in $M^p_m(\rd)$. For the existence
  of a dual window in $M^1_v(\rd)$ see~\cite{GL04}, for the result on Gabor
  expansions see~\cite{book}.  All results of Sections~2 -- 5 carry
  over to general Gabor frames by polarization.
\end{remark}

\section{Further Remarks}

We conclude with an example and some remarks on how Gabor frames are
transformed under FIO.

We first  compute explicitly the multiplier symbol \eqref{symb} for the
dilation operator in dimension $d=1$
$$
D_s f(x)=f(sx) = \int _{\bR } e^{2\pi i s x \eta } \hat{f}(\eta ) \,
d\eta \qquad   s>0 \, .
$$
This is a Fourier integral operator with symbol $\sigma \equiv 1$
and phase $\Phi (x,\eta ) = s x \eta $. Using \eqref{cantra},  the
canonical transformation $\chi$ is calculated to be
$$
\chi(y,\eta )=(y/s,s \eta )\, .
$$
For the lattice  $\Lambda=\a \bZ\times\beta \bZ$ with $\a\beta<1$
and  $\mu=(\a k,\beta \ell)$, $k,\ell \in\bZ$, we then have
$\chi'(\mu)=(\a \lfloor k/s\rfloor,\beta\lfloor s \ell\rfloor)$.
Consequently, the symbols are \begin{align*}{\bf a}_\nu(\mu)&=
\langle T \pi (\mu )g , \pi (\chi  ' (\mu ) +\nu )g\rangle\\& =
e^{2\pi i \a \beta \lfloor \frac{k}s\rfloor\ell'}\int_{\R}e^{2\pi
i \beta( s \ell-\lfloor s\ell\rfloor-\ell')t}g(s t-\alpha
k)\bar{g}(t-\a\lfloor\frac{k}s\rfloor-\a k')\,dt.
\end{align*}
Now, we consider the case of the window function $g(t)=e^{-\pi
t^2}$. A straightforward, but lengthy  computation (as in
~\cite{CF78})  gives
$${\bf a}_{(\a k',\beta \ell')}(\a k,\beta
\ell)=\frac{e^{2\pi i \a \beta (\lfloor
\frac{k}s\rfloor\ell'+(s\ell-\lfloor s \ell\rfloor-\ell')(s
k+\lfloor\frac {k}{s}\rfloor+k'))}}{\sqrt{s^2
+1}}e^{-\pi\beta^2(s\ell-\lfloor sl\rfloor-\ell')^2}e^{-\pi\a^2
s^2(\frac{k}{s}-\lfloor\frac k s \rfloor+k')^2}.
$$

The representation of the dilation operator $D_s$ by a (sum of)
modified Gabor multipliers might be of interest for the numerical
approximation of $D_s$ within the context of \tfa , when one is forced
to use Gabor frames.  Similar formulas  can be worked out
for general metaplectic operators as well.

\vspace{ 3mm}

\emph{More on Gabor frames:} Assume that $\cG (g,\Lambda ) $ is a
frame for $\lrd $ and that $T$ is a FIO with canonical transformation
$\chi $ satisfying the standard conditions. Then the transformed
system $T \cG (g,\Lambda ) =
\{ T\pi (\lambda )g: \lambda \in \Lambda \} $ is a frame, \fif\ $T$ is
invertible. However, if we fix the window $g$ and only warp  the \tf\ space
with $\chi $, then we obtain the set $\cG (g, \chi (\Lambda )) = \{
\pi (\chi (\lambda ))g: \lambda \in \Lambda \}$. It is an interesting
problem to determine when the transformed Gabor system is still a
frame.

If $g(t) = e^{-\pi t^2}$ in dimension $d=1$, then $\cG (g, \chi
(\Lambda ))$ is a frame, \fif\ the lower Beurling density $d^- (\chi
(\Lambda )) >1$. This follows from the characterization of
(non-uniform) Gabor frames by Lyubarskii and Seip~\cite{lyub92,seip92}.

On the other hand, if   $T= D_s$ is   the dilation operator,
$\Lambda = \alpha \bZ \times \beta \bZ $,  and $g$ is a window
with compact support, $\supp g \subseteq [-A,A]$, say, then $\chi
(\Lambda ) = \tfrac{\alpha}{s} \bZ\times s\beta \bZ $ is again a
lattice. Choosing $\tfrac{\alpha }{s} >2A$, then different
translates $g(t-s\inv \alpha k), k\in \bZ $ have disjoint support,
and $\cG (g, \chi (\Lambda ))$ cannot be a frame. Thus the frame
property of $\cG (g, \chi (\Lambda ))$ is subtle.

Based on coorbit theory~\cite{fg89jfa} one can formulate a qualitative
result. Recall that  a (non-uniform) set $\Lambda\subset\R^{2d}$ is
relatively  separated
if $\max _{k\in \zdd } \mathrm{card}\, \Lambda \cap (k+[0,1]^{2d})
<\infty $,  and  $\Lambda $ is called $\delta$-dense for  $\delta>0$, if $
\bigcup_{\lambda\in\Lambda} \overline{B_{\delta}(\lambda)}=\rdd. $

The results in~\cite{fg89jfa} assert that for every $g\in M^1(\rd
)$ there exists a $\delta >0$ depending only on $g$, such that
\emph{every} relatively separated $\delta $-dense set $\Lambda $
generates a frame $\cG (g,\Lambda )$.  We fix $g\in M^1(\rd)$ and
$\delta = \delta (g)>0$.

Now let $T $ be an FIO with canonical transformation $\chi $  with
Lipschitz constant $L$. \emph{If $\Lambda \subseteq \rdd $
 is $\delta /L$-dense, then $\cG (g,\chi (\Lambda ))$ is a frame. }

To see this, observe that
\begin{align*}
\inf _{\lambda \in \Lambda } |z-\chi (\lambda )| & =   \inf _{\lambda
 \in \Lambda } |\chi (\chi \inv  (z)) -\chi (\lambda )| \\
& \leq L \inf _{\lambda \in \Lambda } |\chi \inv (z)-\lambda | \leq L
\, \frac{\delta}{L} = \delta \, .
\end{align*}
Consequently $\chi (\Lambda )$ is $\delta$-dense, and so $\cG (g,
\chi (\Lambda ))$ is a frame.

 \bibliographystyle{abbrv}
 \bibliography{general,new}

\def\cprime{$'$} \def\cprime{$'$} \def\cprime{$'$} \def\cprime{$'$}
  \def\cprime{$'$}

\end{document}